\documentclass[a4paper]{amsart}
\usepackage[T1]{fontenc}
\usepackage{lmodern}
\usepackage[utf8]{inputenc}
\usepackage{amssymb}
\usepackage[all]{xy}
\usepackage{enumitem}
\usepackage{mathrsfs,dsfont,mathtools}
\usepackage{nicefrac}
\usepackage{color}
\usepackage{verbatim}

\usepackage[pdftitle={Equivariant correspondences and the Borel-Bott-Weil theorem},
  pdfauthor={Heath Emerson and Robert Yuncken},
  pdfsubject={Mathematics; MSC 19K35, 46L80}
]{hyperref}

\usepackage[lite]{amsrefs}
\usepackage{microtype}
\numberwithin{equation}{section}
\theoremstyle{plain}
\newtheorem{theorem}[equation]{Theorem}
\newtheorem{lemma}[equation]{Lemma}
\newtheorem{proposition}[equation]{Proposition}
\newtheorem{corollary}[equation]{Corollary}
\theoremstyle{definition}
\newtheorem{definition}[equation]{Definition}

\theoremstyle{remark}
\newtheorem{remark}[equation]{Remark}
\newtheorem{example}[equation]{Example}

\DeclareMathOperator{\KK}{KK}

\DeclareMathOperator{\K}{K}
\newcommand{\Rep}{\mathrm{R}}

\DeclareMathOperator{\RK}{RK}

\newcommand*{\GKK}{\widehat{\textsc{kk}}{}}

\newcommand*{\defeq}{\mathrel{:=}}

\newcommand{\pr}{\mathrm{pr}}
\newcommand{\ID}{\mathrm{id}}

\newcommand{\Ad}{\mathrm{Ad}}

\newcommand{\lie}[1]{\mathfrak{#1}}  
\newcommand{\pt}{\star}

\newcommand{\Spinc}{\mathrm{Spin}^c}

\newcommand{\Weights}{\Gamma_W}
\newcommand{\Roots}{\Delta}
\newcommand{\nvec}{X}
\newcommand{\unitnvec}{Y}

\newcommand{\exteriorC}[1][\bullet]{{\textstyle \bigwedge^{#1}_\C}} 
\newcommand{\length}{l} 
\newcommand{\Hodge}{\circledast}
\newcommand{\Killing}[1]{B(#1)}

\newcommand{\diffeo}{\simeq}
\newcommand{\rep}{\sigma}

\newcommand{\fibration}{\mathrm{p}}
\newcommand{\otherfibration}{\mathrm{q}}
\newcommand{\Thom}[1]{\tau(\fibration_{#1})}
\newcommand{\otherThom}[1]{\tau(\otherfibration_{#1})}
\newcommand{\projection}{\pi}

\newcommand{\zerosection}{\zeta}
\newcommand{\bundleiso}{\varphi}
\newcommand{\otherbundleiso}{\psi}

\newcommand{\cliffrep}{\mathsf{c}}
\newcommand{\Thomendo}{\mathsf{C}}
\newcommand{\otherThomendo}{\mathsf{C}'}

\newcommand{\ptmap}{\mathrm{pt}}

\newcommand{\End}{\mathrm{End}} 

\newcommand*{\CONT}{\mathcal{C}}

\newcommand*{\C}{\mathbb C}
\newcommand*{\Z}{\mathbb Z}

\newcommand*{\N}{\mathbb N}

\newcommand{\Repring}[1]{\mathrm{R}(#1)} 
\newcommand*{\Tvert}{\textup{T}} 

\newcommand*{\mapr}{f}
\newcommand*{\epi}{\twoheadrightarrow}

\newcommand*{\BBW}{\Lambda}

\newcommand{\twisting}[1]{[[#1]]}

\newcommand*{\total}[1]{\lvert#1\rvert}

\newcommand{\Index}{\mathrm{Index}}

\newcommand{\half}{{\frac{1}{2}}}

\newcommand{\ip}[1]{\langle#1\rangle}

\begin{document}

\title{Equivariant correspondences and the Borel-Bott-Weil Theorem}

\author{Heath Emerson}
\email{hemerson@math.uvic.ca}
\address{Department of Mathematics and Statistics\\
  University of Victoria\\
  PO BOX 3045 STN CSC\\
  Victoria, B.C.\\
  Canada V8W 3P4}

\author{Robert Yuncken}
\email{yuncken@math.univ-bpclermont.fr}
\address{Université Blaise Pascal, Laboratoire de Mathématiques\\
  Campus des Cézeaux\\
  B.P. 80026\\
  63171 Aubière cedex\\
  France}

\thanks{Heath Emerson was supported by a Natural Science and 
Engineering Council of Canada (NSERC) Discovery grant. Part of this research was carried out while Robert Yuncken was supported by a PIMS post-doctoral fellowship at the University of 
Victoria.}

\maketitle

\begin{abstract}
We prove an analogue of the Borel-Bott-Weil theorem in equivariant 
\(\KK\)-theory by constructing certain canonical equivariant correspondences 
between minimal flag varieties \(G/B\), with \(G\) a complex semisimple Lie group.

 \end{abstract}

\section{Introduction}

Let $G$ be a complex semisimple Lie group and $B\subset G$ a minimal parabolic subgroup.  Let $\mu$ be a weight for $G$ and $E_\mu$ the corresponding induced holomorphic line bundle on the holomorphic manifold $G/B$.  The Dolbeault cohomology group \(H^*(G/B, E_\mu)\) with its canonical action of \(G\), is a graded-finite-dimensional representation of \(G\), and, more relevantly for us, of its maximal compact subgroup \(K\subset G\). The 
{\em Borel-Bott-Weil} Theorem computes this representation. In this article, we prove a close relative of the Borel-Bott-Weil theorem in the context of equivariant \(\K\)-theory. The bridge between Dolbeault cohomology and \(\K\)-theory is provided by index theory of elliptic operators: 
\(H^*(G/K, E_\mu)\), as a virtual \(K\)-representation, is the same as the \(K\)-index of the Dolbeault operator, $\Index_K [\bar\partial]_\mu \in \Rep (K)$ in the sense of Atiyah and Singer \cite{Atiyah-Singer:Index_I} twisted by 
\(E_\mu\).  From the point of view of Kasparov theory, \([\bar\partial]_\mu\) is a class in the \(K\)-equivariant \(\K\)-homology \(\KK^K(G/B, \C)\) of \(G/B\). This \(\Rep (K)\)-module is acted on by the bivariant group \(\KK^K(G/B, G/B)\), for which a topological model was developed in \cite{Emerson-Meyer:Geometric_KK}, defined using the theory of equivariant \emph{correspondences.} 
The correspondence theory is the main tool used in this work. We construct certain canonical correspondences \(\Lambda (w)\), parameterized by the elements of the Weyl group \(W\), compute how these correspondences act on equivariant \(\K\)-homology, and deduce part of the Borel-Bott-Weil theorem.

The Borel-Bott-Weil Theorem can be broken down into two parts.  Firstly, one can compute $H^*(G/B,E_\mu)$ in the case where $\mu$ is a dominant weight, showing that \(H^i(G/B, E_\mu)\) vanishes in dimensions 
\(i>0\) while in dimension \(i=0\) it is the (unique) irreducible representation of the maximal compact subgroup \(K\) with highest weight \(\mu\). This degree $0$ calculation is effectively the classical Borel-Weil Theorem.  As a second part, one amplifies this calculation by providing a 
degree-altering isomorphism between the cohomology groups $H^*(G/B,E_\mu)$ and $H^*(G/B,E_{\mu'})$, where $\mu'$ is any weight in the orbit of $\mu$ under the shifted (affine) Weyl group action on the weight lattice, explained below. It is this degree-shifting isomorphism which can 
be detected by \(\K\)-theory, as we now explain.

Let $\lie{h}$ be a Cartan subalgebra of the Lie algebra $\lie{g}$ of $G$ and 
$\Weights\subset\lie{h}^*\) the lattice of weights, \emph{i.e.} the group of nonzero 
group homomorphisms \(\lie{h}\to \C^*\). 
 Let $\Roots^+\subset \Weights$ be a choice of a set of positive roots for $G$, which brings with it a generating set of {\em simple reflections} for the Weyl group $W$
 and a corresponding word length function $\length:W\to\N$ on the group. 
 Up to conjugacy, the minimal parabolic subgroup $B\subset G$ is the subgroup with Lie algebra
$\lie{b} = \lie{h} \oplus \bigoplus_{\alpha\in\Roots^+} \lie{g}_\alpha$,
where $\lie{g}_\alpha$ is the $\alpha$-root space of the Lie algebra of $G$.  

Let \([G/B]_\mu\in \KK^K(C(G/B), \C) =: \K_0^K(G/B)\) be the class of the Dolbeault operator 
on \(G/B\) twisted by the \(\K\)-equivariant line bundle \(E_\mu\) obtained by inducing a weight \(\mu\). 
Let $\rho \defeq \half\sum_{\alpha\in\Roots^+} \alpha$ be half the sum of the positive roots.

 \begin{theorem}
 \label{thm:intro:main_theorem}
In the above notation, for any weight $\mu$ of $G$ and any $w\in W$, the identity 
$$ 
 \Index_K [G/B]_\mu = (-1)^{l(w)} \Index_K [G/B]_{w(\mu+\rho)-\rho}  
$$
holds in \(\Rep (K) = KK^K(\C,\C) \). 
 \end{theorem}
  
 The Theorem \ref{thm:intro:main_theorem} is thus a \(\K\)-theoretic version of the 
 Serre duality part of the Borel-Bott-Weil Theorem: a degree-shifting isomorphism 
 between two 
Dolbeault cohomology groups (representations) has been re-stated as an equality, up to a 
specified sign, of (virtual) \(\K\)-representations.

The focus of this article is the $KK$-theory which lies behind Theorem \ref{thm:intro:main_theorem}.
The proof of the theorem
uses the theory of equivariant correspondences and runs as follows.

The Weyl group element $w$ conjugates the subgroup $B$ to another minimal parabolic subgroup $B_w$.  The homogeneous space $G/B\cap B_w$ admits a pair of natural $K$-equivariant holomorphic 
fibrations to $G/B$ and $G/B_w$.  Since the latter space is $K$-equivariantly biholomorphic to $G/B$, we have realized $G/B\cap B_w$ as a holomorphic 
fibred space over $G/B$ in two different ways.  In fact, in each 
case $G/B\cap B_w$ is $K$-equivariantly biholomorphic to the total space of a complex 
vector bundle over $G/B$.  Using the Thom class $\otherThom{w}$ associated to the latter of these 
fibrations we get a 
\(K\)-equivariant holomorphic correspondence 
$$
 G/B \xleftarrow{\fibration_w} (G/B\cap B_w, \otherThom{w}) \xrightarrow{\otherfibration_w} G/B
$$
from \(G/B\) to itself. This yields an element of 
\(\GKK^K(G/B, G/B)\) which we denote by 
\(\BBW (w)\) and call the \emph{Borel-Bott-Weil morphism with 
parameter \(w\in W\).} 

The main result of the paper is the following (they easily imply 
the analogue Theorem \ref{thm:intro:main_theorem} of the Borel-Bott-Weil theorem 
above).

\begin{theorem}
\label{thm:intro:product_formula}
(Borel-Bott-Weil product formula). 
For any weight $\mu$ and $w\in W$, the identity
\begin{equation}
\label{eq:bbw_product}
  \BBW (w)\otimes_{G/B} [G/B]_\mu = (-1)^{\length(w)} [G/B]_{w(\mu + \rho) - \rho}\in 
\KK^K(G/B, \star),
\end{equation}
holds, where \(\BBW (w)\) is the Borel-Bott-Weil 
morphism with parameter $w$.

\end{theorem}

\medskip
We close by noting that we can 
replace \(K\)-equivariance by \(G\)-equivariance 
in Theorem \ref{thm:intro:main_theorem}, using the Baum-Connes conjecture. 
The Borel-Bott-Weil 
theorem is a statement about non-unitary representations of non-compact 
groups. Kasparov theory does not admit such representations. Instead, 
equivariant Kasparov theory for non-compact groups uses unitary, but 
possibly infinite-dimensional representations, and almost-equivariant Fredholm 
operators; these are the cycles for 
the Kasparov representation ring \(\KK^G(\C, \C)\).  There is a restriction map 
\[ \KK^G(A, B) \to \KK^K(A,B)\]
when \(K \subset G\) is a maximal compact subgroup as above, by 
forgetting \(G\)-equivariance to \(K\)-equivariance on cycles. 
The Baum-Connes apparatus shows this map is an isomorphism when \(A\) has the 
form \(A = \CONT(G/B)\otimes A'\) for some \(G\)-C*-algebra \(A'\); this 
follows from a theorem of Tu \cite{Tu:BC_moyennable}. Since all the 
analytic Kasparov classes defined by us have this form, Theorems 
\ref{thm:intro:main_theorem} and \ref{thm:intro:product_formula} have 
their counterparts with \(K\) replaced by \(G\). 

\bigskip




\section{Preliminaries}

\subsection{Equivariant correspondences}
The environment in which the calculations of this paper will 
take place is the topological model for equivariant Kasparov theory 
developed in \cite{Emerson-Meyer:Geometric_KK}. We refer the reader to 
this article for details on the framework. All correspondences used in this 
paper will be smooth, which simplifies the definitions. Let \(K\) be a compact group and let 
\(X\) and \(Y\) be smooth \(K\)-manifolds, \emph{i.e.} smooth manifolds 
with smooth actions of \(K\). 
A smooth correspondence is given by a 
quadruple \( (M, f, b, \xi)\) where 
\begin{itemize}
\item \(M\) is a smooth \(K\)-manifold, 
\item \(f\colon M \to Y\) is a smooth \(K\)-equivariantly \(\K\)-oriented map, 
\item \(b\colon M \to X\) is a smooth \(K\)-equivariant map, and 
\item \(\xi \in \RK^*_{K, X}(M)\) is a smooth \(K\)-equivariant \(\K\)-theory class 
with compact support along the fibres of \(b\) (in the terminology of 
\cite{Emerson-Meyer:Geometric_KK}, a \(\K\)-theory class with \(M\)-compact 
support). 
\end{itemize}

 We usually use the notation
$$
  X \xleftarrow{b} (M,\xi) \xrightarrow{f} Y,
$$
as in \cite{Connes-Skandalis:Longitudinal} (the origin of the theory) to 
denote the quadruple above.

Note that if \(X\) is compact (the case throughout in this article), then 
\(\RK^*_{K, X}(M) = \K^*_K(M)\), the ordinary, compactly supported, 
\(K\)-equivariant \(\K\)-theory of \(M\).

The \emph{degree} of the correspondence is 
the sum of the degrees of \(\xi\) and of \(\mapr\). 

Equivalence classes of equivariant correspondences make up the 
morphisms in the additive category \(\GKK^K\) explained in 
\cite{Emerson-Meyer:Geometric_KK}; it is naturally isomorphic to the corresponding
(analytic) Kasparov ring \(\KK^K\) with some limitations on the arguments; 
these limitations are not relevant for any of the \(K\)-spaces that appear in this 
paper. We generally operate in \(\GKK^K\) in this paper.

 For any pair of 
\(K\)-spaces \(X\) and \(Y\), \(\GKK^K_*(X, Y)\) will thus denote the 
abelian group of equivalence classes of equivariant correspondences 
from \(X\) to \(Y\), graded by parity of degree.  

Two standard examples of $\GKK^K$-classes are important; to fix notion, we 
recall them. 

\begin{example}
\label{ex:right_way_correspondence}
If $b:Y\to X$ is a {\em proper} \(K\)-equivariant map, we define
$$
  b^* \;\defeq\;\;  \big[ X \xleftarrow{b} (Y,\mathbf{1}_Y) \xrightarrow{\ID} Y \big],
$$
where $\mathbf{1}_Y$ is the class of the trivial line bundle $Y\times\C$, the unit in $\RK^*_{K,X}(Y) = \RK^*_K(Y)$.
\end{example}

\begin{example}
\label{ex:wrong_way_correspondence}
If $\Phi$ is an equivariantly $\K$-oriented smooth map from $X$ to $Y$, where \(X\) and \(Y\) are smooth \(K\)-manifolds, we define the {\em wrong-way class} of $\Phi$ as
$$
  \Phi_! \;\defeq\;  \big[ X \xleftarrow{\ID} (X,\mathbf{1}_X) \xrightarrow{\Phi} Y \big],
$$
where $\mathbf{1}_X$ is the class of the trivial line bundle $E\times\C$ in $\RK^*_{K}(X)$.  \end{example}

By a {\em complex $K$-manifold} we shall mean a smooth complex manifold $X$ equipped with a holomorphic action of $K$. The tangent bundle \(\Tvert X\) has a canonical \(K\)-equivariant complex structure and a corresponding \(K\)-equivariant \(\K\)-orientation. This 
supplies an equivariant \(\K\)-orientation on the map from \(X\) to a point.
The corresponding 
wrong-way class is called the \emph{(topological) fundamental class of \(X\)}, and denoted by $[X]$.
Its image in \(\KK^K_0(C_0(X), \C)\) is the class of the \(\K\)-equivariant Dolbeault operator on 
\(X\). 

\medskip

Next, let \(M_1, M_2, Y\) be complex $K$-manifolds. Assume that both \(M_1\) and \(M_2\) are normal $K$-manifolds and that \(\Tvert Y\) is subtrivial. 

Two smooth maps \(f_1\colon M_1\to Y\) and \(b_2\colon M_2 \to Y\) are
\emph{transverse} if for every pair of points \(m_1\in M_1\) and 
\(m_2\in M_2\) with $f_1(m_1)=b_2(m_2)$, the map \(T_{m_1}M_1\oplus T_{m_2}M_2\to 
T_{f_1(m_1)}Y\), \( (\xi_1, \xi_2) \mapsto D_{m_1}f_1 (\xi_1) 
+ D_{m_2}b_2 (\xi_2)\) is surjective.  It is shown in \cite{Emerson-Meyer:Geometric_KK} that when transversality holds, the fibre product 
\[ M_1\times_Y M_2 \defeq \{ (m_1, m_2) \mid f_1(m_1) = b_2(m_2)\}\]
is itself a smooth manifold and the projection \(\pr_2\colon M_1\times_Y M_2 \to M_2\) 
admits a canonical equivariant \(\K\)-orientation from the 
\(\K\)-orientation on \(f_1\). 

If $f_1$ and $b_1$ are {\em holomorphic} maps, the fibre product $M_1\times_Y M_2$ will be a complex manifold, and the projection $\pr_2$ will be a holomorphic map; the 
corresponding \(\K\)-orientation agrees with the one described in the previous 
paragraph. 

\subsection{Complex semisimple Lie groups}
\label{subsec:complex_semisimple_lie_groups}

Here we review some standard structure theory for semisimple groups and fix notation for the remainder of the paper.  For details, see, for example, \cite{Knapp:Lie_groups_beyond_an_introduction}.   

Let $G$ be a complex connected semisimple Lie group and $\lie{g}$ its Lie algebra.  Denote by $\Killing$ its Killing form.  Let $\theta$ be a Cartan involution on $\lie{g}$, so that
$$
  \ip{v,w} \defeq -\Killing{\theta(v),w}, \qquad v,w\in\lie{g}
$$
is a positive definite inner product on $\lie{g}$;  the archetypal example is the operation of negative-conjugate-transpose on $\mathfrak{sl}_n(\C)$.  The $+1$-eigenspace of $\theta$ is the Lie algebra $\lie{k}$ of a maximal compact subgroup $K$ of $G$.

Fix $\lie{h}$, a $\theta$-stable Cartan subalgebra.  Let $\lie{t} = \lie{h}\cap\lie{k}$, which is the Lie algebra of a maximal torus $T$ in $K$.  We have $\lie{h} = \lie{t} \oplus \lie{a}$, where $\lie{a} = i\lie{t}$, and we let $A$ denote the subgroup of $G$ with Lie algebra $\lie{a}$.

The set of roots will be denoted $\Roots$, with $\lie{g}_\alpha$ denoting the root space of $\alpha\in\Roots$.
We fix a choice of positive roots $\Roots^+$, and recall that every positive root is a non-negative integral combination of simple roots.  The lattice of weights will be denoted $\Weights$, and the dominant weights are those $\lambda \in \Weights$ for which $\ip{\lambda,\alpha} \geq 0$ for every positive root $\alpha$.
We will frequently abuse notation by blurring the distinction between a weight $\mu\in\Weights$, the corresponding representation of $T$, and the corresponding holomorphic representation of $H=TA$.

The Weyl group is $W = N_G(H)/Z_G(H)$.  We will frequently identify elements $w\in W$ with a lift to an element of $N_G(H) \subseteq G$, at least when the choice of lift makes no difference.  The usual action of the Weyl group on weights will be denoted by $\mu \mapsto w(\mu)$.  Let $\rho\defeq\half\sum_{\alpha\in\Roots^+}\alpha$ be the half-sum of positive roots.  We will often refer to the {\em shifted action} of the Weyl group, which is the action:
\begin{equation}
\label{eq:shifted_Weyl_action}
  w : \lambda \mapsto w(\lambda+\rho) - \rho.
\end{equation}

We fix the standard Borel subalgebra $\lie{b} \defeq \lie{h} \oplus \lie{n}$, where $\lie{n}$ is the nilpotent subalgebra $\lie{n} \defeq \bigoplus_{\alpha\in\Roots^+} \lie{g}_{\alpha}$.  The associated subgroups are denoted $B$ and $N$.  For each element $w$ of the Weyl group, there are conjugate subgroups
$$
  B_w \defeq wBw^{-1}, \qquad N_w \defeq wNw^{-1}
$$
with corresponding Lie algebras $\lie{b}_w$ and $\lie{n}_w$.  We also define the Lie algebra $\bar{\lie{n}} \defeq \theta\lie{n} = \bigoplus_{\alpha\in\Roots^+} \lie{g}_{-\alpha}$, as well as its conjugates $\bar{\lie{n}}_w \defeq \Ad[w)\bar{\lie{n}}$ for each $w\in W$.

The flag variety of $G$ is the complex homogeneous space $G/B$.  It is $K$-equivariantly diffeomorphic to $K/T$ via the map
$$
  \iota : K/T \xrightarrow{\cong} G/B; \qquad kT \mapsto kB.
$$
However, we shall try to distinguish the two spaces as much as possible.  The difference is technical but important: $G/B$, having a natural complex structure, is canonically $\K$-oriented, while $K/T$ only inherits a $\K$-orientation once it is identified with $G/B$.  Moreover, for any $w\in W$, there is a $K$-equivariant diffeomorphism
$$
  \iota_{w} : K/T \xrightarrow{\cong} G/B_{w}; \qquad kT \mapsto kB_{w},
$$
each inducing a different $\K$-orientation on $K/T$.  This technicality is of course absolutely central to what follows.


\section{The Borel-Bott-Weil theorem}

\subsection{Twisted Fundamental Classes}

Let $\mu$ be a weight of $G$.  As mentioned above, it corresponds to a holomorphic representation of $H$, and one can extend it to a holomorphic character of $B$ which is trivial on $N$.  We denote the one-dimensional representation space by $\C_{\mu}$.  

We shall use the notation $E_{\mu}$ throughout to denote the induced $G$-equivariant line bundle
$$
  E_{\mu} \defeq G\times_{B}\C_{\mu}.
$$
We also have $E_{\mu} \cong K\times_{T}\C_{\mu}$ by restriction.

Recall (see, eg, A.~Wassermann's Frobenius Reciprocity Theorem \cite{Blackadar}[Theorem 20.5.5]) that $\K^*_K(K/T)$ is isomorphic to $\K^*_T(\C) = R(T)$, the representation ring of $T$, as a $\Z$-module.  The representation ring is just $\Z[\Weights]$, the group ring of the weight lattice, and the isomorphism is given by induction:
  
\begin{eqnarray*}
  \mathrm{Ind}_{T}^{K}: ~R(T) &\xrightarrow{\cong} & \K^*_K(K/T)\\
  ~[\mu] &\mapsto& [E_\mu]
\end{eqnarray*}

\begin{definition}
\label{def:twisting}
  Given $\mu\in\Weights$, we define the {\em $\mu$-twisting class} to be the element $\twisting{\mu} \in \GKK^K(G/B,G/B)$ given by the following correspondence:
\begin{equation}
 \label{eq:twisting}
  G/B \xleftarrow{\ID} (G/B, [E_\mu]) \xrightarrow{\ID} G/B.
\end{equation}

  The {\em $\mu$-twisted fundamental class of $G/B$}, denoted $[G/B]_\mu$, is the class of the $K$-equivariant correspondence
\begin{equation}
 \label{eq:twisted_fundamental_class}
  G/B \xleftarrow{\ID} (G/B, [E_\mu]) \rightarrow \star
\end{equation}
 in $\GKK^K(G/B,\star)$.
\end{definition}

Thus, $[G/B]_{\mu} = \twisting{\mu} \otimes_{G/B} [G/B]$, where $[G/B] \defeq [G/B]_{0}$ is the (untwisted) fundamental class of $G/B$.  The reason for the terminology is that $[G/B]$ institutes a duality isomorphism (see \cite{Emerson-Meyer:Geometric_KK}) 
\[ \GKK^K_*(G/B\times X, Y) \cong \GKK^K_*(X, G/B\times Y)\]
valid for arbitrary \(K\)-spaces \(X\) and \(Y\).   
For example if \(X = Y = \star\) then duality gives an isomorphism 
\[ \Repring{T} = \Z[\Gamma_W]
 \cong \GKK^K(G/B, \star).\] This duality can easily be verified 
to send point mass at a weight \(\mu \in \Gamma_W\) to 
the class \([G/B]_\mu\).

\subsection{Borel-Bott-Weil Correspondences}
\label{sec:BBW_correspondences}

Let $w$ be an element of the Weyl group $W=N_G(H)/Z_G(H)$.  
Recall that the subgroup $B_w\defeq wBw^{-1}$ is independent of the choice of lift of $w$ to $N_{G}(H) \subseteq G$.  It is another minimal parabolic subgroup of $G$.   

Consider the homogeneous space $G/(B\cap B_w)$.  This admits two $G$-equivariant fibrations, given by the natural maps 
\begin{eqnarray*}
  p_w: (G/B\cap B_w) &\longrightarrow& G/B, \\
  \otherfibration_w: (G/B\cap B_w) &\longrightarrow& G/B_w .
\end{eqnarray*}
Viewing $G/(B\cap B_w)$ as a $K$-space by restriction, both of these fibrations can be realized as $K$-equivariant vector bundle projections, as we now describe.

Recall that we define $N_w \defeq wNw^{-1}$, $\bar{N}_w \defeq  w\bar{N}w^{-1}$, with Lie algebras $\lie{n}_w$ and $\bar{\lie{n}}_w$ respectively.  Then $\lie{n} =  (\lie{n}\cap\bar{\lie{n}}_w) \oplus (\lie{n}\cap\lie{n}_w) $ is a decomposition of $\lie{n}$ into Lie subalgebras.  Since $N$ is a connected simply-connected nilpotent Lie group, there is a corresponding factorization $N = (N\cap \bar{N}_w)(N\cap N_w) $.

\begin{lemma}
\label{lem:bundle_isomorphism}
Let $w \in W$.  One can define a $K$-equivariant diffeomorphism\linebreak $\bundleiso_w:\total{K\times_T (\lie{n}\cap\bar{\lie{n}}_w)}\xrightarrow{\cong} G/(B\cap B_w)$ by the formula $\bundleiso_w:[k,\nvec] \mapsto k\exp(\nvec).(B\cap B_w)$ such that the 
diagram  
\begin{equation}
\label{eq:bundle_isomorphism1}
\xymatrix{
   K\times_T (\lie{n}\cap\bar{\lie{n}}_w)   \ar[r]_-{\cong}^-{\bundleiso_w}  \ar[d]_{\projection_w}  & G/(B \cap B_w) \ar[d]_{\fibration_w} \\
  K/T \ar[r]_{\cong} & G/B
}
\end{equation}
commutes.  Moreover, $\bundleiso_w$ is fibrewise holomorphic (with respect to the fibrations $\projection_{w}$ and $\fibration_{w}$).
\end{lemma}

In other words, the \(K\)-equivariant fibrations 
\( K\times_T (\lie{n}\cap\bar{\lie{n}}_w) \to K/T\) and\linebreak \( G/(B\cap B_{w}) \to G/B\) are 
equivalent in the category of \(K\)-equivariant fibrations with holomorphic fibres.

\begin{proof}

To see that the map $\bundleiso_w$ is well-defined we compute, for any $k\in K$, $t\in T$, $\nvec\in\lie{n}\cap\bar{\lie{n}}_w$:
$$
  \bundleiso_w( [kt,\Ad (t^{-1}) \nvec] ) = kt.t^{-1} \exp (\nvec) t  (B\cap B_w) = k\exp(\nvec).(B\cap B_w) = \bundleiso_w( [k,\nvec] ).
$$

Next, surjectivity.  Let $g\in G$ be arbitrary.  There is a decomposition $G = KNA = K(N\cap \bar{N}_w)(N\cap N_w) A$, and we decompose $g$ as $g=kn_1n_2a$ accordingly.  Since $(N\cap N_w) A \subseteq B\cap B_w$, we have $\bundleiso_w([k, \log(n_1)]) = g(B\cap B_w)$.  

Next suppose $[k,\nvec]$ and $[k',\nvec']\in K\times_M (\lie{n}\cap\bar{\lie{n}}_w)$ have the same image under $\bundleiso_w$.  Since $B\cap B_w = T(N\cap N_w)A$, there exist $t\in T$ and $n_2a\in (N \cap N_w)A$ such that
$$
  k'\exp(\nvec') =k\exp(\nvec)tn_2a = kt\exp(\Ad(t^{-1} \nvec))n_2a.
$$
By the uniqueness of the $K(N\cap \bar{N}_w)(N\cap N_w)A$-decomposition, $k'=kt$ and $\nvec'=\Ad(t^{-1})\nvec$, which is to say $[k',\nvec']=[k,\nvec]$.  

That the map is $K$-equivariant is straightforward, as is the commutativity of the diagram of bundle maps.  Fibrewise holomorphicity follows from the holomorphicity of the exponential map.
\end{proof}

\begin{remark}
\label{rmk:bundle_isomorphism}
There is an alternative realization of the space $G/(B\cap B_w)$ as a $K$-equivariant vector bundle via the diagram
\begin{equation}
\label{eq:bundle_isomorphism2}
\xymatrix{
  K\times_T (\lie{n}_w\cap\bar{\lie{n}})   \ar[r]_-{\cong}^-{\bundleiso'_w}  \ar[d]_{\projection'_w}  & G/(B \cap B_w) \ar[d]_{\otherfibration_w} \\
  K/T \ar[r]_{\cong} & G/B_w
}
\end{equation}
where the top map has essentially the same defining formula: $\bundleiso'_w:[k,\nvec'] \mapsto k\exp(\nvec').(B\cap B_w)$.  The proof is basically identical. Thus,  the holomorphic manifold 
$G/(B\cap B_w)$ admits two distinct structures as a complex \(K\)-equivariant vector 
bundle over \(K/T\), via the maps 
\(\pi_w\) and \(\pi_w'\). This point will be of crucial importance later.
\end{remark}

\begin{definition}
 \label{def:zero_section}
Using the diagrams \eqref{eq:bundle_isomorphism1} and \eqref{eq:bundle_isomorphism2}, 
we may consider the zero sections of the two complex vector bundles 
\(K\times_T (\lie{n}_w\cap\bar{\lie{n}})\) and \(K\times_T (\lie{n}\cap\bar{\lie{n}}_w) \) 
as \(K\)-equivariant 
maps $\zerosection_w:G/B \to G/(B\cap B_w)$ and $\zerosection'_w:G/B_w \to G/(B\cap B_w)$.  They are given simply by 
\begin{eqnarray*}
 \zerosection_w &:& kB \mapsto k(B\cap B_w) \qquad\text{for $k\in K$}, \\
 \zerosection'_w &:& kB_w \mapsto k(B\cap B_w) \qquad\text{for $k\in K$}, 
\end{eqnarray*}
where we stress that in applying these formulas, we are obliged to choose coset 
representatives \(k\) belonging to the compact subgroup $K$. 
\end{definition}

The importance of realizing $G/(B \cap B_{w})$ as a complex $K$-vector bundle over $G/B$ is that there is a Thom class
$$
  \Thom{w} \in \K_K^*(G/(B\cap B_w)),
$$
obtained by pushing forward the Thom class from $\K^*(\total{K\times_T (\lie{n}\cap\bar{\lie{n}}_w)} )$.

Note that this Thom class is dependent upon the fibration map $\fibration:G/(B\cap B_w) \to G/B$. The alternative fibration $\otherfibration_w:G/(B\cap B_w) \to G/B_w$ defines a different class $\otherThom{w}$, pushed forward from $\K_K^* (\total{K\times_T (\lie{n}_w\cap\bar{\lie{n}})} )$.

\medskip

Let $w\in W$.  The spaces $G/B_w$ and $G/B$ are $G$-equivariantly diffeomorphic, even biholomorphic, via the right multiplication map $R_w:g.(wBw^{-1}) \mapsto gw.B$.   We can now define one of our main objects of study.

\begin{definition}
\label{def:BBW_correspondence}
The {\em Borel-Bott-Weil morphism \(\BBW(w) \in \GKK^K(G/B, G/B)\) with parameter 
\(w\in W\)} is the class 
of the \(\K\)-equivariant holomorphic correspondence
\begin{equation}
\label{eq:uninflated_BBW_correspondence}
  \xymatrix{
     G/B &  (G/B\cap B_w,\otherThom{w}) \ar[l]_-{\fibration_w} \ar[r]^-{\otherfibration_w} & G/B_w \ar[r]^{R_w}_\diffeo & G/B.
  }
\end{equation}
\end{definition}

\begin{example}
If $w=e$ is the identity element, then $B\cap B_e=B$, $\lie{n}\cap\bar{\lie{n}}_w$ is the zero Lie subalgebra, inducing to the zero vector bundle on \(K/T\), and $\Thom{e}$ is the Thom class \([1]\)
of the zero vector bundle. Thus $\BBW(e)=1\) is represented by the correspondence 
\[ \xymatrix{ 
G/B & (G/B, [1])\ar[l]^{\mathrm{id}}\ar[r]^{\mathrm{id}} & G/B}\]
which is the identity correspondence. Thus
\( \BBW(e) = 1 \in \GKK^K(G/B, G/B)$.
\end{example}

\subsection{Product structure}

For $w\in W$, we denote by $\iota_w : K/T \xrightarrow{\diffeo} G/B_w$ the \(K\)-equivariant diffeomorphism defined by $kT\mapsto kB_w$ for $k\in K$.

Each of these maps identifies $K/T$ with a complex manifold, with \(K\) acting by a holomorphic action, and this complex structure induces a corresponding
 \(K\)-equivariant $\Spinc$-structure on \(K/T\). 
  All of these $\Spinc$-structures will be different.  To keep track of them, we use the complex picture whenever possible.

\begin{definition}
 We denote by $I_w:G/B \to G/B_w$ the (non-holomorphic) \(\K\)-equivariant 
 diffeomorphism defined by the commuting diagram
$$
 \xymatrix{
   G/B \ar[r]^{I_w} & G/B_w \\
   K/T \ar[u]^{\iota_e}_\diffeo \ar[r]^{\ID} & K/T \ar[u]^{\iota_w}_\diffeo.
 }
$$
Thus, $I_w$ corresponds to the identity map on $K/T$ but with an unusual 
$\K$-orientation.
\end{definition}

If \(w\in W\), then right translation \(R_w\colon G/B_w\to G/B\) is a \(K\)-equivariant 
map yielding an element \(R_w^*\in \GKK^K(G/B, G/B_w)\). 
The following proposition asserts, roughly, that after twisting \(R_w\) 
by the change of equivariant \(\K\)-orientation induced by \(I_w\), we 
get exactly the Borel-Bott-Weil correspondence \(\BBW(w)\).

\begin{proposition}
\label{prop:compact_picture}
 Let $w\in W$.  Then $\BBW(w) = (I_{w}^{-1} \circ R_{w}^{-1})^*$ in $\GKK^K(G/B,G/B)$.
\end{proposition}

\begin{proof}
The map $R_w$ is biholomorphic, so $(R_{w}^{-1})^* = {R_w}_!$. Using the realization of $G/(B \cap B_w)$ as a $K$-equivariant vector bundle over $G/B_w$, we can perform a Thom modification to get
\begin{eqnarray*}
 (I_{w}^{-1} \circ R_{w}^{-1})^* &=& \big[ G/B \xleftarrow{I_{w}^{-1}} G/B_w \xrightarrow{R_w} G/B \big]\\
  &=& \big[ G/B \xleftarrow{I_{w}^{-1} \circ \otherfibration_w} (G/(B\cap B_w), \otherThom{w}) \xrightarrow{R_w\circ\otherfibration_w} G/B \big].
\end{eqnarray*}

Next consider the commutative diagram
$$
 \xymatrix{
  G/B \ar[r]^{I_w} \ar[d]_{\zerosection_w} & G/B_w \ar[d]^{\zerosection'_w} \\
  G/(B\cap B_w) \ar@{=}[r] & G/(B\cap B_w).
 }
$$
Since the bundle projections $\fibration_w$ and $\otherfibration_w$ are homotopy inverses to the zero sections $\zerosection_w$ and $\zerosection'_w$, we have that $I_w^{-1} \circ \otherfibration_w \sim \fibration_w$.  Hence $(I_{w}^{-1} \circ R_{w}^{-1})^* = \BBW(w)$.
\end{proof}

\begin{corollary}
 \label{cor:group_embedding}
 The map $w \mapsto \BBW(w)$ is a group homomorphism from the Weyl group into the invertible elements of the ring $\GKK^K(G/B,G/B)$.
\end{corollary}

\begin{proof}
 One just needs to check that $R_{w_1}\circ I_{w_1} \circ R_{w_2} \circ I_{w_2} = R_{w_1w_2}\circ I_{w_1w_2}$.  This is immediate if one represents elements of $G/B$ as $kB$ with $k\in K$.
\end{proof}

\subsection{Commutation relations in $\widehat{\textsc{\small KK}}\mathbf{{}^K(G/B,G/B)}$}
\label{sec:commutation}

We begin with some generalities on pullbacks of induced bundles. 

Let $H_2 \leq H_1 \leq G$ be a nested sequence of closed Lie subgroups, and let $V$ be a vector space with a representation of $H_1$.  If $\fibration : G/H_2 \epi G/H_1$ denotes the canonical fibration map, then there is an equivariant bundle isomorphism
\begin{equation}
 \label{eq:pullback_of_induced_bundle}
\fibration^* (G\times_{H_1} V) \cong G\times_{H_2} V, 
\end{equation}
given by following pullback diagram:
$$
\xymatrix{
  G\times_{H_2} V \ar[r] \ar[d] & G\times_{H_1} V \ar[d] &\quad&   
      {[g,v]} \ar@{|->}[r] \ar@{|->}[d] &  {[g,v]} \ar@{|->}[d] \\
  G/H_2 \ar[r]^p & G/H_1 &\quad& 
      gH_2 \ar@{|->}[r] & gH_1.
}
$$

\medskip

Recall that each weight $\mu$ defines a one-dimensional holomorphic representation of $B$.  It will be convenient to use an explicit notation for this in the next few paragraphs, so we denote it by $\rep_\mu : B \to \End(\C_\mu)$.  We shall denote by $\rep^w_\mu$ the representation of $B_w$ defined by conjugating by $w\in W$: $\rep^w_\mu(wbw^{-1}) \defeq \rep_\mu(b)$.

Then there is a $G$-equivariant bundle isomorphism 
\begin{equation}
 \label{eq:Weyl_action_on_bundles}
R_w^*(G\times_B \C_\mu) = G\times_{B_w} \C_{\mu},
\end{equation}
where the representation of $B_w$ on the right-hand side is $\rep^w_{\mu}$.
The appropriate pullback diagram is:
$$
\xymatrix{
  G\times_{B_w} \C_{w(\mu)} \ar[r] \ar[d] & G\times_{B} \C_{\mu} \ar[d] &\quad&   
      {[g,v]} \ar@{|->}[r] \ar@{|->}[d] &  {[gw,v]} \ar@{|->}[d] \\
  G/B_w \ar[r]^{R_w} & G/B &\quad& 
      gB_w \ar@{|->}[r] & gwB.
}
$$

\begin{lemma}
 \label{lem:pullback_bundles}
For any $\mu\in\Weights$ and $w\in W$ we have $\otherfibration_w^* R_w^* E_\mu \cong \fibration_w^* E_{w(\mu)}$ as $G$-equivariant complex line bundles over $G/(B\cap B_w)$.
\end{lemma}

\begin{proof}
 As described above, $\otherfibration_w^* R_w^* E_\mu \cong \otherfibration_w^* (G\times_{B_w} \C_{\mu})$.  Restricting the conjugated representation $\rep_\mu^w$ to $B \cap B_w$ yields a representation which is trivial on $N\cap N_w$ and given by $e^{w(\mu)}$ on $TA$.  Thus \eqref{eq:pullback_of_induced_bundle} gives $ \otherfibration_w^* R_w^* E_\mu \cong G\times_{B \cap B_w} \C_{w(\mu)}$.  This is isomorphic to $\fibration_w^* E_{w(\mu)}$ by \eqref{eq:pullback_of_induced_bundle} again.
\end{proof}

\begin{proposition}
 \label{prop:commutation_relation}
For any $\mu\in\Weights$ and $w\in W$, $\BBW(w) \otimes_{G/B} \twisting{\mu} = \twisting{w(\mu)} \otimes_{G/B} \BBW(w)$.
\end{proposition}

\begin{proof}
 We calculate,
\begin{eqnarray*}
 \lefteqn{\BBW(w)\otimes_{G/B}\twisting{\mu}} \\
  \quad &=& G/B \xleftarrow{\fibration_w} (G/B\cap B_w, \otherThom{w}) \xrightarrow{R_w \circ \otherfibration_w} G/B
   \xleftarrow{\ID} (G/B,[E_\mu]) \xrightarrow{\ID} G/B \\
  &=& G/B \xleftarrow{\fibration_w} (G/B\cap B_w, \otherThom{w}.\otherfibration_w^*R_w^*[E_\mu] ) \xrightarrow{R_w \circ \otherfibration_w} G/B \\
  &=& G/B \xleftarrow{\fibration_w} (G/B\cap B_w, \otherThom{w}.\fibration_w^*[E_{w(\mu)}] ) \xrightarrow{R_w \circ \otherfibration_w} G/B \\
  &=& G/B \xleftarrow{\ID} (G/B,[E_{w(\mu]}]) \xrightarrow{\ID} G/B \xleftarrow{\fibration_w} (G/B\cap B_w, \otherThom{w}) \xrightarrow{R_w \circ \otherfibration_w} G/B \\
  &=& \twisting{w(\mu)} \otimes_{G/B} \BBW(w).
\end{eqnarray*}
\end{proof}

\subsection{Comparing Thom classes}
\label{sec:comparing_Thom_classes}

We begin this section by comparing the two Thom classes $\Thom{w}$ and $\otherThom{w}$ on the space $G/B \cap B_w$ (see Section \ref{sec:BBW_correspondences}).  It will suffice to consider the case where $w$ is the reflection in a simple root $\alpha$.  in that case we have $\lie{n}\cap\bar{\lie{n}}_w = \lie{g}_\alpha$ and $\lie{n}_w\cap\bar{\lie{n}} = \lie{g}_{-\alpha}$.  

Recall that $\Thom{w}$ is the pushforward of the Thom class of $\total{K\times_T (\lie{n}\cap\bar{\lie{n}}_w)} = \total{K\times_T \lie{g}_{\alpha}}$ via the bundle isomorphism of \eqref{eq:bundle_isomorphism1}.  Taking advantage of the complex structure on the fibres, the corresponding spinor bundle is $K\times_M \exteriorC \lie{g}_{\alpha}$.  There is an $\Ad(T)$-invariant inner product on $\lie{g}_{\alpha}$ via the Killing form.  Letting $\lambda_\nvec$ denote the exterior product by $\nvec\in\lie{g}_{\alpha}$, we have a Clifford algebra representation
\begin{equation}
\label{eq:Clifford_representation}
  \cliffrep:\lie{g}_{\alpha} \to \End(\exteriorC \lie{g}_{\alpha}) ; \quad \cliffrep(\nvec)\defeq \lambda_\nvec - \lambda_\nvec^*.
\end{equation}
The Thom class of $\total{K\times_T \lie{g}_{\alpha}}$ is the pullback of the spinor bundle  along the bundle projection $\projection_w: K\times_T \lie{g}_\alpha \to K/T$, equipped with the bundle endomorphism which at each point is the Clifford representation of that point.  

Since $K\times_T \exteriorC \lie{g}_{\alpha} \cong \C_0 \oplus \C_{\alpha}$, we can identify the spinor bundle over $K/T$ with $G\times_B (\C_0 \oplus \C_{\alpha})$.  The space $\C_{\alpha}$ here identifies naturally with $\lie{g}_{\alpha}$ as a $T$-space, but not as a $B$-space: we have made an arbitrary extension to a $B$-representation.  

Using Equation \eqref{eq:pullback_of_induced_bundle}, the push-forward of the Thom class by $\otherbundleiso$ is then
\begin{equation}
 \label{eq:Thom_class}
  \Thom{w} = ( G \!\! \underset{B \cap B_w}{\times} \!\!(\C_0 \oplus \C_{\alpha}) , \; 
    \Thomendo_w ),
\end{equation}
where $\Thomendo_w$ is the bundle endomorphism defined at each point of $G/(B \cap B_w)$ by
$$
 \Thomendo_w(k\exp(\nvec)(B\cap B_w)) = \cliffrep(\nvec) \qquad \text{for } k\in K, ~\nvec \in \lie{g}_{\alpha}.
$$

A similar calculation shows that the Thom class $\otherThom{w}$ associated to the other projection is
\begin{equation}
 \label{eq:other_Thom_class}
  \otherThom{w} = ( G \times_{B \cap B_w} (\C_0 \oplus \C_{-\alpha}) , \; 
    \otherThomendo_w ),
\end{equation}
where 
$$
 \otherThomendo_w(k\exp(\nvec')(B\cap B_w)) = \cliffrep(\nvec') \qquad \text{for } k\in K, ~\nvec' \in \lie{g}_{-\alpha}.
$$

\medskip
To compare these two classes, we define a homotopy.  For $t\in[0,1]$, define a map
\begin{eqnarray*}
 \gamma_t : G/(B\cap B_w) &\to& G/(B\cap B_w), \\
   k\exp(\nvec) (B\cap B_w) &\mapsto& k\exp(t\nvec) (B\cap B_w), \qquad \text{for }k\in K, ~\nvec \in \lie{g}_{\alpha}.
\end{eqnarray*}
This is just the pushforward by $\bundleiso_w$ of the retraction of the bundle $K\times_T \lie{g}_{\alpha}$ to the zero section.

Consider the smooth family $\Phi_t$ of bundle endomorphisms of $G \times_{B \cap B_w} (\C_0 \oplus \C_{-\alpha})$ defined by
$$
  \Phi_t(z) \defeq \otherThomendo(\gamma_t(x)), \qquad (z\in G/ B\cap B_w).
$$
Since $\gamma_0$ has image the zero section, $\Phi_0$ is the zero endomorphism.  By smoothness and the compactness of $G/B$, the family
$$
  \Psi_t \defeq \frac{1}{t} \Phi_t  \qquad (t\neq0)
$$
has a well-defined limit at $t=0$, which we denote by $\Psi_0$.  

\begin{lemma}
\label{lem:homotopy}
Let $\theta$ denote the Cartan involution on $\lie{g}$.  At a point\linebreak  $k\exp(\nvec)(B\cap B_w)$ of $G/(B\cap B_w)$, where $k\in K$ and $\nvec \in \lie{g}_{\alpha}$, the limit\linebreak $\Psi_0(k\exp(\nvec)B\cap B_w)$
is the endomorphism of the fibre $\exteriorC \lie{g}_{-\alpha}$ defined by
 $$
  \Psi_0(k\exp(\nvec)B\cap B_w) = \cliffrep(-\theta\nvec) .
 $$ 
\end{lemma}

\begin{proof}
We have $\gamma_t(k\exp(\nvec)(B \cap B_w)) = k\exp(t\nvec)(B \cap B_w)$.  By the Campbell-Baker-Hausdorff formula, $\exp(t\nvec) = \exp(t(\nvec+\theta\nvec))\exp(-t\theta\nvec)\exp(o(t))$.  Since $\exp(t(\nvec+\theta\nvec)) \in K$, we have that $\Psi_t$ acts on the fibre at $k\exp(\nvec)B\cap B_w$ by
$$
 \Psi_t(k\exp(\nvec)B\cap B_w) = \frac{1}{t}\cliffrep(-t\theta\nvec + o(t)).
$$
which has limit $\cliffrep(-\theta\nvec)$ as $t\to0$.
\end{proof}

In the next lemma, we fix identifications of $\lie{g}_{\pm\alpha}$ with $\C$ by identifying some arbitrary unit vector $Y\in\lie{g}_{-\alpha}$ with $1$, and likewise with $\theta Y \in \lie{g}_{\alpha}$.  Ultimately the choice of this $Y$ makes no difference.

\begin{lemma}
Fix $\unitnvec \in \lie{g}_{-\alpha}$ with $\|\unitnvec\| = 1$.  Define a grading-reversing map $\beta : \exteriorC \lie{g}_{-\alpha} \to  \exteriorC \lie{g}_{\alpha}$ by
$$
\beta :
\begin{cases}
  \omega \mapsto \omega.\theta\unitnvec,& \qquad\text{ for } \omega\in\exteriorC^0\lie{g}_{-\alpha} = \C,\\
  \nvec' \mapsto \ip{\unitnvec,\nvec'}, & \qquad\text{ for } \nvec'\in\exteriorC^1\lie{g}_{-\alpha} = \lie{g}_{-\alpha}.
\end{cases}
$$
Then for any $\nvec\in\lie{g}_\alpha$,
$$
  \beta^{-1} \cliffrep(\nvec) \beta = \cliffrep(-\theta\nvec).
$$
\end{lemma}

\begin{remark}
 Equivalently, $\beta = \theta\circ\Hodge$, where $\Hodge$ is the (anti-linear) Hodge $*$-operator on $\exteriorC \lie{g}_{-\alpha}$.
\end{remark}

\begin{proof}
 We calculate
\begin{equation*}
 \begin{array}{rll}
 \beta^{-1}\lambda_\nvec\beta : & \omega \mapsto 0,  & \text{for $\omega\in\C$},\\
 \beta^{-1}\lambda_\nvec\beta : & \omega\unitnvec \xmapsto{\beta} \omega \xmapsto{\lambda_\nvec} \omega\nvec \xmapsto{\beta^{-1}} \ip{\theta\unitnvec,\omega\nvec}, \qquad & \text{for $\omega\unitnvec\in\lie{g}_{-\alpha}$}
 \end{array}
\end{equation*}
and
\begin{equation*}
 \begin{array}{rll}
 \lambda_{\theta X}^* : & \omega \mapsto 0, & \text{for $\omega\in\C$}, \\
 \lambda_{\theta X}^* : & \omega\unitnvec \mapsto \ip{\theta\nvec,\omega\unitnvec}, \hspace{3cm} & \text{for $\omega\unitnvec\in\lie{g}_{-\alpha}$}.
 \end{array}
\end{equation*}
These maps are equal 
since $\theta$ is anti-unitary. 
Also $\beta^{-1}\lambda_\nvec^*\beta = \lambda_{\theta\nvec}$, by the unitarity of $\beta$.  The result now follows from the definition $\cliffrep(\nvec) \defeq \lambda_\nvec - \lambda_\nvec^*$.
\end{proof}

The map $\beta$ is not $T$-equivariant --- it alters the weights, since it maps $\lie{g}_{-\alpha}$ to $\C_0$ and $\C_0$ to $\lie{g}_\alpha$.  But if we alter it by defining
\begin{eqnarray*}
   \beta' : \exteriorC \lie{g}_{-\alpha} &\to&  (\exteriorC \lie{g}_{\alpha}) \otimes \lie{g}_{-\alpha} \\
     Z &\mapsto& \beta Z \otimes \unitnvec,
\end{eqnarray*}
then it is weight-preserving, and hence $T$-equivariant.  It induces a grading-reversing bundle isomorphism
\begin{multline}
  \ID \times_{B\cap B_w} \beta' : G \times_{B\cap B_w} \exteriorC \lie{g}_{-\alpha} \to
    G \times_{B\cap B_w} ((\exteriorC \lie{g}_{\alpha}) \otimes \lie{g}_{-\alpha}) \\
  \cong (G \times_{B\cap B_w} \exteriorC \lie{g}_{\alpha}) \otimes_{G/B\cap B_w} \fibration^*E_{-\alpha},
\end{multline}
which intertwines the bundle endomorphisms $\Psi_0$ and $\Thomendo_{w} \otimes \ID$.  Combining this with the fact that $\Psi_1 = \Thomendo_{w}$, we have proven the following fact.

\begin{proposition}
 \label{prop:comparing_Thom_classes}
If $w\in W$ is the reflection in the simple root $\alpha$, then $\otherThom{w} = -\Thom{w}\otimes\fibration^*[E_{-\alpha}]$ in $\K_K^*(G/(B\cap B_w))$.
\end{proposition}

\begin{remark}
 There is a more general formula:  for any $w\in W$, $\otherThom{w} = -\Thom{w}\otimes\fibration^*[E_{w(\rho) -\rho}]$, where $\rho$ is the half-sum of the positive roots.  This can be proven along the same lines as above with significantly more work, or deduced from results to follow.  We shall not need it.
\end{remark}

\subsection{The Borel-Bott-Weil Theorem: Action on $\K$-homology and Indices}
\label{sec:K-homology}

\begin{proof}[Proof of Theorem \ref{thm:intro:product_formula}]
We wish to show
$$
  \BBW(w) \otimes_{G/B} [G/B]_\mu = (-1)^{\length(w)}\,[G/B]_{w(\mu+\rho)-\rho}.
$$
By the multiplicativity of the map $w\mapsto\BBW(w)$, it suffices to take $w$ a reflection in a simple root $\alpha$.  

Let $\mu\in\Weights$.  Using the fact that $[G/B]_\mu = \twisting{\mu}\otimes_{G/B}[G/B]$, Proposition \ref{prop:commutation_relation} gives
$$
 \BBW(w) \otimes_{G/B} [G/B]_\mu = \twisting{w(\mu)} \otimes_{G/B} \BBW(w) \otimes_{G/B} [G/B] \\
$$
From Proposition \ref{prop:comparing_Thom_classes},
\begin{eqnarray*}
 \BBW(w) \otimes_{G/B} [G/B] &=& \left[ G/B \xleftarrow{\fibration_w} (\, G/(B\cap B_w) , \, -\Thom{w}\fibration^*[E_{-\alpha}] \,) \rightarrow \star \right] \\
  &=& \left[ G/B \xleftarrow{\ID} ( G/B , \, -[E_{-\alpha}] ) \rightarrow \star \right] \\
  &=& -[G/B]_{-\alpha}.
\end{eqnarray*}
So we get
$$
 \BBW(w) \otimes_{G/B} [G/B]_\mu =  -[G/B]_{w(\mu)-\alpha}.
$$
Since $w$ is the reflection in $\alpha$, we have $\alpha = w(\rho)-\rho$, which proves the result.

\end{proof}

We now pass to the index-theoretic application. 
Let $\ptmap:G/B\to\pt$ denote the map of $G/B$ to a point and 
\(\ptmap^*\in \GKK^K(\C, G/B)\) its topological $\KK$-theory class. 

 For a weight $\mu$,  the topological \(\K\)-index 
 of the twisted fundamental class 
 $[G/B]_\mu\in\GKK^G(G/B,\star )$ is defined by
$$
  \Index_K [G/B]_\mu \defeq \ptmap^*\otimes_{G/B} [G/B]_\mu \quad \in \GKK^K(\C,\C).
$$
We do not bother to use different notation for the \emph{analytic} index \(\Index_K [G/B]_\mu \in \KK^K(\C, \C) \cong \Rep (K)\); which one we are talking about will be made clear by the context.  The analytic index is the same as the cohomology group \(H^*(G/B, E_\mu)\) figuring in the classical Borel-Bott-Weil theorem, while the Atiyah-Singer index theorem asserts that it is equal to the image of the topological index under the map \(\GKK^K (G/B, \star) \to \KK^K(\CONT(G/B), \C)\).

\begin{proof}[Proof of Theorem \ref{thm:intro:main_theorem}]
We note that a Thom modification yields
\begin{eqnarray}
\label{eq:index_of_BBW_class}
 \ptmap^* \otimes_{G/B} \BBW(w) 
   &=& \left[ \star \leftarrow (G/(B\cap B_w), \, \otherThom{w}) \xrightarrow{ R_w \circ \otherfibration_w} G/B \right] \\
   &=& \left[ \star \leftarrow G/B_w \xrightarrow{ R_w } G/B \right] \nonumber \\
   &=& \ptmap^*. \nonumber 
\end{eqnarray}
Composing with the $\mu$-twisted fundamental class on the right, and applying Theorem \ref{thm:intro:product_formula} gives $(-1)^{\length(w)} \Index_K [G/B]_{w(\mu+\rho)-\rho} = \Index_K [G/B]_\mu$.
\end{proof}

\begin{remark}
\label{rmk:K-theory_action}

Let us also record the action of the Borel-Bott-Weil classes on equivariant $\K$-theory.  The induction isomorphism $\Repring{T} \xrightarrow{\cong} \K_K(K/T)$ associates to $[\mu]$ the correspondence
$$
  [E_\mu] \defeq \left[ \star \leftarrow (K/T,\, [E_\mu]) \xrightarrow{\ID} K/T \right] = \ptmap^* \twisting{\mu}.
$$
Thus, if we compose the commutation relation of Proposition \ref{prop:commutation_relation} on the left by $\ptmap^*$ and use Equation \eqref{eq:index_of_BBW_class}, we get the right action:
$$
  [E_{w(\mu)}] \otimes_{G/B} \BBW(w) = [E_\mu].
$$

\end{remark}

\bibliographystyle{alpha}

\bibliography{BBW}

\begin{bibdiv}
\begin{biblist}

\bib{Atiyah-Singer:Index_I}{article}{
      author={Atiyah, M.~F.},
      author={Singer, I.~M.},
       title={The index of elliptic operators. {I}},
        date={1968},
        ISSN={0003-486X},
     journal={Ann. of Math. (2)},
      volume={87},
       pages={484\ndash 530},
      review={\MR{0236950 (38 \#5243)}},
}

\bib{Blackadar}{book}{
      author={Blackadar, Bruce},
       title={{$K$}-theory for operator algebras},
     edition={Second},
      series={Mathematical Sciences Research Institute Publications},
   publisher={Cambridge University Press},
     address={Cambridge},
        date={1998},
      volume={5},
        ISBN={0-521-63532-2},
      review={\MR{1656031 (99g:46104)}},
}

\bib{Connes-Skandalis:Longitudinal}{article}{
      author={Connes, A.},
      author={Skandalis, G.},
       title={The longitudinal index theorem for foliations},
        date={1984},
        ISSN={0034-5318},
     journal={Publ. Res. Inst. Math. Sci.},
      volume={20},
      number={6},
       pages={1139\ndash 1183},
         url={http://dx.doi.org/10.2977/prims/1195180375},
      review={\MR{775126 (87h:58209)}},
}

\bib{Emerson-Meyer:Geometric_KK}{article}{
      author={Emerson, Heath},
      author={Meyer, Ralf},
       title={Bivariant {$K$}-theory via correspondences},
        date={2010},
        ISSN={0001-8708},
     journal={Adv. Math.},
      volume={225},
      number={5},
       pages={2883\ndash 2919},
         url={http://dx.doi.org/10.1016/j.aim.2010.04.024},
      review={\MR{2680187 (2012a:19013)}},
}

\bib{Knapp:Lie_groups_beyond_an_introduction}{book}{
      author={Knapp, Anthony~W.},
       title={Lie groups beyond an introduction},
     edition={Second},
      series={Progress in Mathematics},
   publisher={Birkh\"auser Boston Inc.},
     address={Boston, MA},
        date={2002},
      volume={140},
        ISBN={0-8176-4259-5},
      review={\MR{1920389 (2003c:22001)}},
}

\bib{Tu:BC_moyennable}{article}{
      author={Tu, Jean-Louis},
       title={La conjecture de {B}aum-{C}onnes pour les feuilletages
  moyennables},
        date={1999},
        ISSN={0920-3036},
     journal={$K$-Theory},
      volume={17},
      number={3},
       pages={215\ndash 264},
         url={http://dx.doi.org/10.1023/A:1007744304422},
      review={\MR{1703305 (2000g:19004)}},
}

\end{biblist}
\end{bibdiv}

\end{document}